\documentclass[11pt]{amsart}
\usepackage{amssymb}
\usepackage{amsfonts}
\usepackage{amscd}
\usepackage{amsmath}
\usepackage{mathrsfs}
\usepackage{curves}
\usepackage{epsfig}
\usepackage[pass]{geometry}
\usepackage{graphicx}
\usepackage{verbatim}
\usepackage{latexsym}
\usepackage{eucal}

\numberwithin{equation}{section}

\newtheorem{theorem}{Theorem}[section]

\newtheorem{prop}[theorem]{Proposition}
\newtheorem{cor}[theorem]{Corollary}

\newtheorem*{result}{Main Results}

\def \mca {{\mathscr A}}
\def \mcb {{\mathscr B}}
\def \mcc {{\mathscr C}}
\def \mcd {{\mathscr D}}
\def \mce {{\mathscr E}}
\def \mcf {{\mathscr F}}

\def \mci {{\mathscr I}}

\def \mcl {{\mathscr L}}
\def \mcm {{\mathscr M}}

\def \mcp {{\mathscr P}}

\def \mcr {{\mathscr R}}
\def \mcs {{\mathscr S}}

\def \mbr {{\mathbb R}}

\def \id {\operatorname{Id}}
\def \comp {\operatorname{comp}}
\def \loc {\text{loc}}

\def \diag{\textrm{Diag}}
\def \supp {\text{supp }}

\def \beqq {\begin{equation}}
\def \eeqq {\end{equation}}
\def \WF {\text{WF}}

\def \bpf {\begin{proof}}
\def \epf {\end{proof}}
\def \beq {\begin{equation*}}
\def \eeq {\end{equation*}}

\def \eps {\epsilon}   
\def \la {\lambda}   
\def \La {\Lambda}    
   
\def \lap {\Delta}
\def \p {\partial}
\def \ha {\frac{1}{2}}

\def \tilde {\widetilde}
\def \hat {\widehat}
\def \fnf {\frac n4}
\def \diam {\operatorname{diam}}

\begin{document}
\title{Convolutional neural networks in phase space and inverse problems}
\author{Gunther Uhlmann}
\address{Gunther Uhlmann
\newline
\indent Department of Mathematics, University of Washington 
\newline
\indent and Institute for Advanced Study, the Hong Kong University of Science and Technology}
\email{gunther@math.washington.edu}
\author{Yiran Wang}
\address{Yiran Wang
\newline
\indent Department of Mathematics, Stanford University}
\email{yrw@stanford.edu}
\begin{abstract}
We study inverse problems consisting on 
determining medium properties using the responses to probing waves from the machine learning point of view. Based on the understanding of propagation of waves and their nonlinear interactions, we construct a deep convolutional neural network in which the parameters are used to classify and reconstruct the coefficients of nonlinear wave equations that model the medium properties. Furthermore, for given approximation accuracy, we obtain the depth and number of units of the network and their quantitative dependence on the complexity of the medium.  
\end{abstract}

\maketitle

\section{Introduction}\label{sec-intro}
In this work, we consider inverse problems for nonlinear hyperbolic equations. The method to be developed applies to a large class of hyperbolic equations on manifolds, however, for simplicity, we 
consider nonlinear acoustic wave equations on $\mbr^3$ of the form
\beqq\label{eqnon}
\begin{gathered}
(\p_t^2 + c^2(x) \lap ) u(t, x) + F(t, x, u(t, x)) = f(t, x),\quad t > 0, \ \ x\in \mbr^3\\
u(t, x) = 0, \quad t \leq 0, \ \ x\in \mbr^3
\end{gathered}
\eeqq
where $c(x) > 0$ is the wave speed, $f(t, x)$ is the source term, $\lap = -\sum_{i = 1}^3 \p_{x^i}^2$ is the (positive) Laplacian on $\mbr^3$ and $F(t, x, u)$ is a smooth function in $t, x$ and $u$. We are mainly interested in the case when $F$ is nonlinear in $u.$ We denote
\beq
P = \p_t^2  + c^2(x) \lap
\eeq
the linear wave operator. Because we only consider local problems later, we assume that $c(x)$ is non-trapping without loss of generality. One can think of  equation \eqref{eqnon} as modeling acoustic waves generated by the source $f(t, x)$ traveling in a medium with certain nonlinear mechanism. The coefficients $c(x), F(t, x, u)$ characterizes the medium properties. The inverse problem (to be formulated precisely in Section \ref{sec-inv}) we address is the determination of the wave speed $c(x)$ and the nonlinear term $F(t, x, u)$ by measuring the response of waves traveling through the medium. From the machine learning point of view, the problem is to learn material properties (characterized by $c$ and $F$) from the data  (the source and wave responses).  In particular, our goal is to classify different materials from the data but furthermore, we aim to reconstruct the materials from the data. 

The linearized problem, that is when $F(t, x, u)$ is linear in $u$ and typically with the hyperbolic Dirichlet-to-Neumann data, has been studied extensively in the literature. There are well-developed methods such as Boundary Control (BC) method, see \cite{KKM} for an overview. However, the nonlinear problem to be considered in this work are not always solvable by linearization. Some recent progress have been made towards solving these problems by exploiting the nonlinear interactions of waves, beginning with the work by Kurylev-Lassas-Uhlmann \cite{KLU}. The phenomena that nonlinear interactions of waves could generate new waves have been known for a while and observed in many physical experiments. Mathematically this phenomena has been studied from the point of view of interactions of singularities by the notable work of Bony \cite{Bo}, Melrose-Ritter \cite{MR}, Reed-Rauch \cite{RR}, S\'a Barreto \cite{Sa}, Melrose-S\'a Barreto-Zworski \cite{MSZ}, Zworski \cite{Zw} etc. See also Beals \cite{Bea} for an overview of the subject in the 80s and 90s.  The idea introduced in \cite{KLU} is that by using distorted plane waves concentrated near fixed directions, one can keep track of their interactions and the newly generated waves. One can characterize the ``features'' of these waves in the data which eventually leads to the determination the parameters of the equation. 

Inspired by these ideas, we construct a deep neural network for solving the inverse problem of recovering $c$ and $f$ from the measured data 
 and we prove approximation properties of the network. The coefficients $c(x), F(t, x, u)$ can be reconstructed from the parameters of the network. The informal version of our main theorem is 
\begin{result}
We construct a deep convolutional neural network (in Section \ref{sec-conv}, see Figure \ref{figdnn}) with $M$ levels, $K$ units on each level  and parameter set $\Theta$ such that for data $(f, u)$ where $u$ is the solution of \eqref{eqnon} with source $f$ and $\|f\| < \eps$, the network generates an approximation function $h(f; \Theta)$ satisfying 
\[
\|u - h(f; \Theta)\| < C_M \eps^M
\]
The norms are specified in Theorem \ref{thmmain}. The number $K$ and $M$ depends on the complexity of $c(t, x), F(t, x, u)$. The parameters $\Theta$ can be used to reconstruct $c(x), F(t, x, u)$ (see Section \ref{sec-app}). 
\end{result}


Roughly speaking, the units in each level of the neural network represent small wave units. Units in deeper levels capture the effects of wave propagation and nonlinear interactions. In fact, for general source term $f$ (not necessarily distorted plane waves), we think of it as consisting of sufficiently many small wave units and the network captures the interaction among them. We will see that it is natural to work in the phase space and consider the high frequency information in the wave units, which we take as the ``features'' in this machine learning problem. A main part of the construction is to show how these ``features'' interact with each other and propagate through the network. This is a new feature even for the conventional convolutional network and Mallat's scattering network \cite{Ma1}. Another novelty of our result is that we study the depth and number of units of the network in terms of the complexity of the parameter functions $c(x), F(t, x, u)$. Roughly speaking, for highly nonlinear functions, a deeper network should be used to reveal such effects and produce better approximations. Finally, we interpret the meaning of the  parameters in the network and explain how they can be used for reconstructing $c(x), F(t, x, u).$

The method we developed should certainly be applicable to other inverse problems or machine learning problems involving wave equations. We should mention at this point that we focus on theoretical questions about the network properties in this work.  The analysis of the resulting optimization problem and numerical implementation will be pursued elsewhere. 

The organization of the paper is follows. In Section \ref{sec-inv}, we formulate the inverse problem to be considered in this article. Then we compare the iteration scheme for solving wave equations (Section \ref{sec-ite}) and the deep forward networks (Section \ref{sec-mlp}). We propose our first network in Section \ref{sec-net1} and discuss some issues in the architecture related to the wave phenomena. The resolution of these issues leads us to the convolutional neural network in Section \ref{sec-conv} but before that, we need to discuss the nonlinear interactions of conormal waves (Section \ref{sec-con}),  the estimates of linear propagation (Section \ref{sec-lin}) and nonlinear effects (Section \ref{sec-non}). Finally, we prove the approximation properties of the network and reconstruct the coefficients in Section \ref{sec-app}. 

\section{The inverse problem}\label{sec-inv}
We consider two types of inverse problems for equation \eqref{eqnon} with different types of data. In this work, we shall work with the first problem exclusively but we remark that the same methods apply to the second problem as well. 

\subsection{The source perturbation problem}
Let $f(t, x)$ be compactly supported. For $T>0$ fixed, consider \eqref{eqnon}
\beq
\begin{gathered}
(\p_t^2 + c^2(x) \lap ) u(t, x) + F(t, x, u) = f(t, x),\quad t \in (-\infty, T],\quad  x\in \mbr^3\\
u(t, x) = 0, \quad t\leq 0, \quad x\in \mbr^3
\end{gathered}
\eeq
It is known (also see Section \ref{sec-ite}) that for $f \in H^s([0, T]\times \mbr^3)$ sufficiently small and compactly supported, there is a unique solution $u\in  H^{s+1}([0, T]\times \mbr^3)$. We denote this solution map by $u = L(f)$. We remark that we do not pursue the optimal regularity result in this work. 

\begin{figure}[htbp]
\centering
\includegraphics[scale = 0.7]{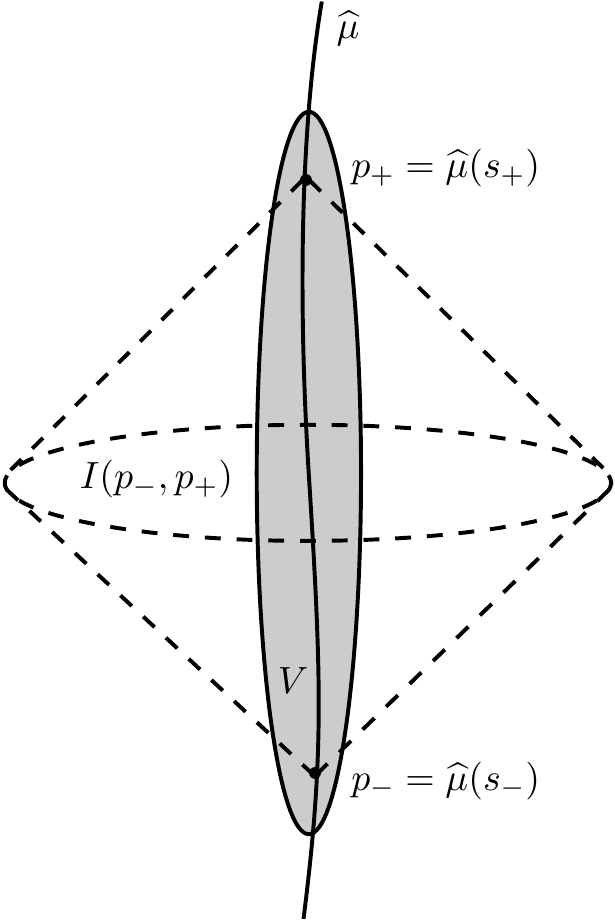} \quad \quad \quad \quad \quad \quad \quad \quad
\includegraphics[scale = 0.7]{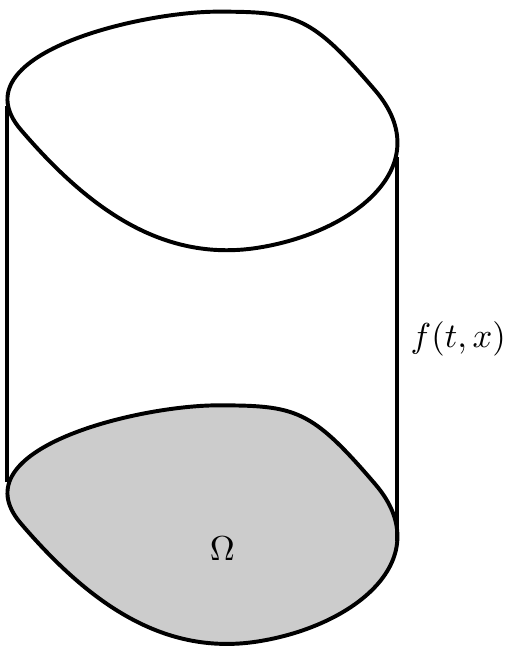}
\caption{Two types of inverse problem. Left: the source perturbation problem. Right: the boundary value problem.}
\label{figinv}
\end{figure}

 We want to determine $c$ and $F$ in the region where the wave can travel to. It is convenient to formulate the problem using the space-time nature of wave propagation. Let 
\[
g = dt^2 + c^{-2}(x)dx^2
\]
 be the Lorentzian metric so that the corresponding Laplace-Beltrami operator is giving by  $\square_g = \p_t^2 + c^2(x) \lap$. We denote $\mcm = \mbr^4$ and consider the Lorentzian manifold $(\mcm, g)$. Let $\hat \mu(s) \subset \mcm$ be a time-like geodesic where $s\in [-1, 1]$. In general relativity, this represents the world line of a freely falling observer. Let $V\subset \mcm$ be an open relatively compact neighborhood of $\hat \mu([s_-, s_+])$ where $-1<s_-<s_+<1$. We denote $\mcm(T) = (-\infty,T]\times \mbr^3$ and choose $T>0$ such that $V\subset \mcm(T)$. Let $p_\pm = \hat \mu(s_\pm)$. See the left of Figure \ref{figinv} .

We recall some notions of causalities, see e.g.\ \cite{BEE}.  For $p, q\in M$, we denote by $p\ll q$ ($p\leq q$) if $p\neq q$ and $p$ can be joined to $q$ by a future pointing time-like (causal) curve. We denote by $p\leq q$ if $p = q$ or $p<q$. The chronological future of $p\in M$ is the set $I^+(p) = \{q\in M: p\ll q\}$. The causal future of $p\in M$ is $J^+(p) = \{q\in M: q\leq p\}$. The chronological past and causal past are denoted by $I^-(p)$ and $J^-(p)$ respectively. For any set $A\subset M$, we denote the causal future by $J^\pm(A) = \cup_{p\in A}J^\pm(p)$. Also, we denote $J(p, q) = J^+(p)\cap J^-(q)$ and $I(p, q) = I^+(p)\cap I^-(q)$. 

Let $f$ be supported in $V$ and we measure the wave $u$ in $V$. The data set is 
\[
\mcd_{sour} \doteq \{(f, u|_V): u = L(f), f\in H^s_{comp}(V), s> 1\}. 
\]
The inverse problem is to determine $c(x)$ and $F(t, x, u)$ on $I(p_-, p_+)$ from $\mcd_{sour}$, see Figure \ref{figinv}. Notice that $I(p_-, p_+)$ is the largest set that the wave $u$ can travel to from $V$ and return to $V$. 

This formulation was introduced for the Einstein equations in \cite{KLU1} which has a concrete physical interpretation, that is to determine space-time structures (e.g.\ topological, differentiable structure and the metric) from actively generated gravitational perturbations measured near a freely falling observer. In fact, the Einstein equation in wave gauge is a second order quasilinear hyperbolic system. The problem has been further studied in \cite{LUW2} for Einstein-Maxwell equations and \cite{UW} for more general source fields. One of our motivations  is to develop an algorithm to understand the gravitational wave interactions in these work. For semilinear wave equations on globally hyperbolic Lorentzian manifolds, the problem was studied in \cite{KLU} and \cite{LUW1}.

\subsection{The hyperbolic Dirichlet-to-Neumann problem.}
For the second type of inverse problem the information is given on the boundary. We consider the wave equation \eqref{eqnon} on a bounded domain $\Omega\subset \mbr^3$ with smooth boundary $\p \Omega$. See the right of Figure \ref{figinv}. For fixed $T>0$, consider 
\beq
\begin{gathered}
(\p_t^2  + c^2(x) \lap ) u(t, x) + F(t, x, u(t, x)) = 0,\quad (t, x)\in [0, T]\times \Omega \\
u(t, x) = f(t, x), \quad  t \leq T, x\in \p \Omega,\\
u(t, x) = 0, \quad t\leq 0, x\in \Omega.
\end{gathered}
\eeq
For $f\in H^s([0, T]\times \p \Omega)$ sufficiently small and regular, and  compactly supported, the problem is well-posed. See for example \cite{DH} for the treatment of Cauchy data and also \cite{DUW, NW}. We can define the Dirichlet-to-Neumann map 
\[
\La(f) = \nu\cdot \p u|_{[0, T]\times \p \Omega} 
\]
where $\nu$ is the outward normal vector to $\p \Omega$. The data set is 
\[
\mcd_{DtN}\doteq \{(f, \La(f)): f\in H_{comp}^s([0, T]\times \p \Omega), s > 1\}
\]
The inverse problem is to determine $c(x)$ and $F(t, x, u)$ from this data set. We remark that on unbounded domain, one can formulate the problem as a scattering problem. 

For this setup, Nakamura-Watanabe \cite{NW} considered the one dimensional quasilinear wave equation, which is further generalized in Nakamura-Vashisth \cite{NV} for systems in one dimension. The nonlinear elastic system is of particular interest because of its applications in geophysics and rock sciences. For example, one is interested in determining the underground formation of the Earth using nonlinear responses of seismic waves because the contrast in nonlinear parameters are stronger than linear ones, see \cite{KSC, TR1}.  In de Hoop-Uhlmann-Wang \cite{DUW}, the authors analyzed the nonlinear interaction of two elastic waves and the inverse problems of determining elastic parameters is addressed as well. 

\section{The iteration scheme}\label{sec-ite}
We begin with the iteration method for solving nonlinear wave equations. This material is rather classical, however, we want to show to what extent each iteration step reveals nonlinear effects. To illustrate the idea, we take the polynomial nonlinear function
\[
F(t, x, u) = a(t, x) u^2 + b(t, x) u^3 + c(t, x) u^4
\]
as an example. The coefficients $a, b, c$ reflects the nonlinearity in increasing orders.  Also, we shall consider the small source perturbation problem for the wave equation
 \beqq\label{eqnonlin}
 P u(t, x) + F(t, x, u) = \eps f(t, x), \quad  (t, x)\in \mcm(T),
 \eeqq
 where $f$ is compactly supported and $\eps$ is a small parameter. These two simplifications will be removed eventually. 
 
Let $v$ be the solution of the linearized equation on $\mcm(T)$
\[
Pv = f
\]
It is well-known that there is a fundamental solution $Q = P^{-1}$. We write $v = Q(f)$. Let $u$ be the solution of \eqref{eqnonlin}. Then formally we have 
\beq
\begin{gathered}
P(u - \eps v) + F(u) = 0 
\Longrightarrow u  = \eps v - Q(F(u)). 
\end{gathered}
\eeq
Here, we omitted the dependence of $F$ on $t, x$ in the notation. Now we let $u^{(1)} = \eps v$ be the linearized solution and set 
\beq
\begin{split}
u^{(2)}  &= \eps v - Q(F(u^{(1)}))\\
 &= \eps v - \eps^2 Q(av^2) + O(\eps^3).
\end{split}
\eeq
 We observe that modulo $O(\eps^3)$ terms, the  coefficients  $a$ appear in $u^{(2)}$ and this is associated with the quadratic nonlinearity. We continue this procedure to get 
\beq
\begin{split}
u^{(3)}  &= \eps v - Q(F(u^{(2)}))\\
 &= \eps v - \eps^2 Q(av^2)  + 2 \eps^3 Q(avQ(av^2)) -  \eps^3  Q(bv^3)  + O(\eps^4)
\end{split}
\eeq
and another iteration gives 
\beq
\begin{split}
u^{(4)}  &= \eps v - Q(F(u^{(3)}))\\
 &= \eps v - \eps^2 Q(av^2) + 2 \eps^3 Q(avQ(av^2))  -  \eps^3  Q(bv^3) + \\
  & \eps^4 [- Q(cv^4) + 2 Q(avQ(bv^3))  +  3 Q(bv^2Q(av^2)) -  4 Q(avQ(avQ(av^2)))] + O(\eps^5).
\end{split}
\eeq
The point is that for each $i = 1, 2, 3 $, modulo $O(\eps^i)$ terms, we should expect to see the nonlinear coefficients in $u^{(i)}$.  
One continue the procedure to obtain that the sequence $u^{(n)}$. The fact is that $u^{(n)}$ converges to the solution $u$ in a proper sense.  
 
 \begin{prop}\label{propite}
 Consider the nonlinear wave equation 
 \beq
 \begin{gathered}
Pu(t, x) + F(t, x, u(t, x)) = f(t, x),\quad (t, x) \in \mcm(T)\\
u(t, x) = 0, \quad (t, x)\in \mcm(0).
\end{gathered}
\eeq
We assume that $F$ is a smooth function with $F(t, x, 0) = F_u(t, x, 0) = 0$. Fixed $T> 0$, there exists $\eps_0$ such that for $f$ compactly supported in $\mcm(T)\backslash \mcm(0)$ with $\| f\|_{H^s(\mcm)}\leq \eps, s>1, 0< \eps < \eps_0$, the sequence $u^{(n)}$ defined iteratively by 
\[
u^{(1)} = Q(f), \quad u^{(n)}  = u^{(0)} - Q(F(t, x, u^{(n-1)})), n \geq 2
\]
converges to a unique solution $u\in H^{s+1}(\mcm(T))$. Moreover,  we have the estimates 
\[
\|u^{(n)} - u\|_{H^{s+1}(\mcm(T))} < C_n \eps^n, \quad \|u\|_{H^{s+1}(\mcm(T))} < C\eps
\]
where $C_n, C>0$ depends on $c, F$ and $C_n$ depends on $n$ as well.  
 \end{prop}
 \bpf
First of all, we recall that $Q: H^{s}_{comp}(\mcm(T))\rightarrow H^{s+1}_{loc}(\mcm(T))$ is bounded, see for example \cite[Prop.\ 5.6]{DUV}. So there is $C_Q > 0$ depending on $c$ such that  
\[
\|Qf\|_{H^{s+1}}\leq C_Q\|f\|_{H^s}. 
\]
For $s > 1$, the space $H^{s+1}(\mcm(T))$ is an algebra. Moreover, $F(t, x, u)\in H^{s+1}(\mcm(T))$ for any smooth function $F$ and $u\in H^{s+1}(\mcm(T))$, see \cite{Tay}.  By Sobolev embedding, $H^{s+1}(\mcm(T))\subset C^{r}(\mcm(T))$ with $r < s - 1$. In particular, $u^{(n)} \in H^{s+1}(\mcm(T))\subset C^0(\mcm(T))$ are continuous for $s > 1.$

We want to show that $u^{(n)}$ form a Cauchy sequence. For convenience, we take $u^{(0)} = 0.$ Suppose $f$ is supported in a compact set $K\subset \mcm(T)$. By finite speed of propagation for linear wave equations, we know that each $u^{(n)}, n \geq 1$ is supported in $J_+(K)$. We shall assume $u^{(n)}$ supported in $J_+(K)\cap \mcm(T)$. 

Now we consider $u^{(m)} - u^{(n)}, m, n \geq 1$ satisfying 
 \beq
 \begin{gathered}
P(u^{(m)} - u^{(n)})= -[ F(t, x, u^{(m-1)}) - F(t, x, u^{(n-1)})],\quad (t, x)\in \mcm(T)\\
u^{(m)} - u^{(n)} = 0, \quad (t, x) \in \mcm(0).
\end{gathered}
\eeq
We obtain 
\beqq\label{eqest1}
\begin{gathered}
\|u^{(m)} - u^{(n)}\|_{H^{s+1}} \leq C_Q \|F(t, x, u^{(m-1)}) - F(t, x, u^{(n-1)})\|_{H^{s}}.
\end{gathered}
\eeqq

First we take $n = 1$ to get 
\beq 
\begin{gathered}
\|u^{(m)} - u^{(1)}\|_{H^{s+1}} \leq C_Q \|F(t, x, u^{(m-1)})\|_{H^{s}}.
\end{gathered}
\eeq 
Then we write 
\beq
F(t, x, u^{(m-1)}) = (\ha \int_0^1 \p_u^2F(t, x, \tau u^{(m-1)}) dt) (u^{(m-1)})^2. 
\eeq
Because $F$ is smooth, and $u^{(n)}\in H^{s+1}$, we can use Moser type estimates (see \cite[Prop.\ 3.9]{Tay2} which also works for $F(t, x, u)$ by minor modifications of the proof), to obtain that for $(t, x)\in J_+(K)$ and $\tau \in [0, 1]$, 
\beq
\| \p_u^2F(t, x, u^{(m-1)}) - \p_u^2 F(t, x, 0)\|_\infty \leq C \|u^{(m-1)}\|_\infty (1 + \|u^{(m-1)}\|_{H^{s+1}})
\eeq
where $C$ depends on $|\p_u^kF(t, x, u)|$ for $k\leq s+1$.  
 Thus,
\beq
\| \p_u^2F(t, x, u^{(m-1)})\|_\infty \leq C_F + C \|u^{(m-1)}\|_\infty (1 + \|u^{(m-1)}\|_{H^{s+1}})
\eeq
and we have 
\beq 
\begin{gathered}
\|u^{(m)} - u^{(1)}\|_{H^{s+1}} \leq C_Q[ C_F + C \|u^{(m-1)}\|_\infty (1 + \|u^{(m-1)}\|_{H^{s+1}})]\| u^{(m-1)} \|^2_{H^{s+1}} 
\end{gathered}
\eeq 
Now we use induction and assume that $\|u^{(m-1)} - u^{(1)}\|_{H^{s+1}} < \eps $ for $\eps $ sufficiently small. This implies that $\|u^{(m-1)}\|_{H^{s+1}}\leq C_0\eps$ for some constant $C_0.$ We see that 
\beqq\label{eqest2}
\|u^{(m)}- u^{(1)}\|_{H^{s+1}} \leq \eps \bigg( \eps  C_0^2  C_Q[ C_F + C \eps  (1 +\eps )] \bigg)
\eeqq
So we just need to take $\eps< \eps_0$ with $\eps_0 C_0^2C_Q[ C_F + C \eps_0  (1 +\eps_0 )] < 1$ and we obtain $\|u^{(m)} - u^{(1)}\|_{H^{s+1}} < \eps $. This finishes the induction and  shows that 
\beqq\label{equm}
\|u^{(m)}\|_{H^{s+1}}\leq C_0 \eps
\eeqq
are bounded for all $m.$

Next, we return to \eqref{eqest1} and write
\beq
 F(t, x, u^{(m-1)}) - F(t, x, u^{(n-1)}) = (\int_0^1 \p_u F(t, x, u^{(m-1)} + \tau u^{(n-1)}) dt) (u^{(m-1)} - u^{(n-1)})
\eeq
As $\p_uF(t, x, 0) = 0$, we use Moser type estimate again to get 
\beq
\begin{gathered}
\|\p_u F(t, x, u^{(m-1)} + \tau u^{(n-1)})\|_\infty \leq C\|u^{(m-1)} + \tau u^{(n-1)}\|_\infty (1 + \|u^{(m-1)} + \tau u^{(n-1)}\|_{H^{s+1}}) 
\leq C_1 \eps 
\end{gathered}
\eeq
for all $\tau \in [0, 1]$. 
Now we use the fact that $u^{(n)}$ are continuous and bounded on $J_+(K)$ to get 
\beq
\begin{gathered}
\|u^{(m)} - u^{(n)}\|_{H^{s+1}} \leq C_QC_FC_1\eps \|u^{(m-1)} - u^{(n-1)}\|_{H^{s+1}} 
\end{gathered}
\eeq
which implies that for $m>n$
\beq
\|u^{(m)} - u^{(n)}\|_{H^{s+1}} \leq (C_QC_FC_1 \eps)^{n-1} \|u^{(m-n+1)} - u^{(1)}\|_{H^{s+1}} \leq (C_QC_FC_1 \eps)^{n-1} \eps
\eeq
where we have used \eqref{eqest2}. By possibly shrinking $\eps_0$ further so that $C_QC_FC_1\eps < 1$ for $\eps< \eps_0$, we see that  $u^{(m)}$ is a Cauchy sequence and it converges to some $u\in H^{s+1}$. Then we have the estimates  
\beq
\|u^{(n)} - u\|_{H^{s+1}} \leq  C_n \eps^{n}
\eeq
for some constant $C_n$ depending on $n$. The estimates of $\|u\|_{H^{s+1}}$ follows from triangle inequality and \eqref{equm}.
 \epf
 
We make a few remarks. (1) The same argument works for quasilinear wave equations. (2) The proof works for short time (i.e.\ for small $T$) instead of for small data. (3) Here we mainly consider $F$ nonlinear. If $F$ has linear terms, that is $F_u(t, x, 0)\neq 0$, then the same arguments work through. We just have to replace $P$ by $\tilde P = P + F_u(t, x, 0)u$ and change $Q$ to $\tilde Q = \tilde P^{-1}$. In case $F(t, x, 0)\neq 0$, one can solve $Pw = -F(t, x, 0)$ first and repeat the argument to get 
\[
\|u^{(n)} - (u - w)\|_{H^{s+1}(\mcm(T))}\leq C_n \eps^n. 
\]
So in principle, one can remove the small data assumption. For simplicity, we stick with this assumption in the rest of the paper.

\section{Deep feedforward networks}\label{sec-mlp}
The iteration scheme for solving wave equations shares some similarities to the architecture in deep feedforward networks. We briefly review it and refer the reader to \cite{Go} for details.

In general, the goal of deep feedforward networks is to approximate some function $y = F^*(\vec x): \mbr^N\rightarrow \mbr^M$. For example, in classification problems, the function returns the number of classes the data $\vec x$ belongs to. The feedforward network defines a mapping $y = F(\vec x; \vec \theta)$ where $\vec\theta \in \mbr^M$ is the parameter, and learns the value of the parameter that gives an approximation of the function $F^*$.

There are many variants of deep forward networks. We illustrate using the multi-layer perceptrons (MLPs). The construction of MLP consists of a sequence of composition of linear mappings (the perceptron) followed by nonlinear maps called the activation function. Usually, the linear mapping is taken to be the affine transformation
\[
f(\vec x;  \theta) = A \vec x + \vec b
\]
where $A\in \mbr^N\times \mbr^N$ and $\vec b\in \mbr^N$ are the parameters $\theta = \{A, \vec b\}$.  There are many choices of the activation function in practice. A commonly used one is the rectified linear unit (ReLU) 
\[
g(z) = \max\{0, z\}, \quad z\in \mbr.
\]

To build the MLP, we start from the input data $\vec x$ and call it $h^{(0)} = \vec x$. This forms the first layer of the network. 
Let $\theta^{(1)} = \{A^{(1)}, \vec b^{(1)}\}$ be the first set of parameters. We set
\[
h^{(1)} = g(f(h^{(0)}; \theta^{(1)})) = g( A^{(1)} h^{(0)} + \vec b^{(1)})
\]
where $g$ applies to each components of $f(h^{(0)};  \theta^{(1)})$. This defines the second layer, also called the first hidden layer of the network. It is worth mentioning that one can introduce multiple affine transformations on the same level, that is 
\[
h^{(1)} = g( \sum_{k = 1}^K \tilde h_k^{(1)}), \quad \tilde h_k^{(1)} = A_k^{(1)} h^{(0)} + \vec b_k^{(1)}
\]
and each $\tilde h_k^{(1)}$ is called a unit for this level.  One must realize that without the activation function, $h^{(1)}$ would be just a linear function of $h^{(0)}$. 

This defines the iteration scheme. We continue to obtain 
\[
h^{(n)} = g(f(h^{(n-1)};  \theta^{(n)})), \quad n = 1, 2, \cdots, N, 
\]
which defines the $n$-th layer of the network. Here, $N$ is called the depth of the network. Eventually, we obtain the approximation function   $F(\vec x; \vec \theta) = h^{(N)}$ where $\vec \theta$ is the collection of parameters $\theta^{(n)}, n = 1, 2, \cdots, N.$ 

To find the parameters $\vec\theta$, we need to solve an optimization problem on some training data set $X$. The cost function can be formulated as
\[
J(\theta) = \sum_{\vec x\in X} ||F^*(\vec x) - F(\vec x; \vec \theta)||^2
\]
in the $l^2$ norm for vector spaces.  The cost function is usually non-convex and the problem is solved by gradient descent based method such as stochastic gradient descent and back-propagation, see \cite{Go, LeC}.

On the theoretical level, it is important to understand the approximation ability of the network. It is shown in \cite{Ho1, Ho2, Le} that MLPs can approximate any Borel measurable function, which is known as the universal approximation property. We recall and formulate two such theorems below with special attention to the regularity of activation functions. The first theorem is for continuous activation functions and the second one is for $C^l$ activation functions. 

Let $G: \mbr\rightarrow \mbr$ be a Borel measurable function. Let $\mca$ be the set of affine transformations, namely $A\in \mca$ means $Ax = w  x + b, x\in \mbr$. We define
\[
\Sigma(G) = \{f: \mbr^n \rightarrow \mbr | f(x) = \sum_{j =1}^q \beta_j G(A_j(x)), x\in \mbr, \beta_j\in \mbr, A_j\in \mca, q = 1, 2, \cdots.\}
\]

\begin{theorem}[Theorem 2.1 of \cite{Ho1}]\label{thmuni1}
Let $G: \mbr\rightarrow \mbr$ be a continuous non-constant function.  Then for any compact set $K\subset \mbr$, $\Sigma(G)$  is dense in $C^r$ with respect to 
\[
\rho_K(f, g) = \sup_{x\in K}|f(x) - g(x)|
\] 
\end{theorem}

\begin{theorem}[Corollary 3.4 of \cite{Ho2}]\label{thmuni2}
Let $G: \mbr\rightarrow \mbr$ be a function in $C^l(\mbr^n)$ with non-negative integer $l$ satisfying 
\beqq\label{eqcond}
\int |\frac{d^l}{dx^l}G| dx < \infty.
\eeqq
 Then for any compact set $K\subset \mbr$, $\Sigma(G)$  is dense in $H^{m}(\mbr^n)$ for $m\leq l$ with respect to 
\[
\rho^m_K(f, g) = \sum_{j = 0}^m \sup_{x\in K}|\frac{d^j}{dx^j}f(x) - \frac{d^j}{dx^j}g(x)|
\] 
\end{theorem}

Another important variant of feedforward networks is the convolutional neural network where  the affine transformation is replaced by convolutions. See \cite{LeC, Ma2}. We refer readers to \cite{Go} for the motivation and advantages of this type of network. 

\section{A network for solving inverse problems}\label{sec-net1}
\begin{figure}[htbp]
\centering
\includegraphics[scale=.8]{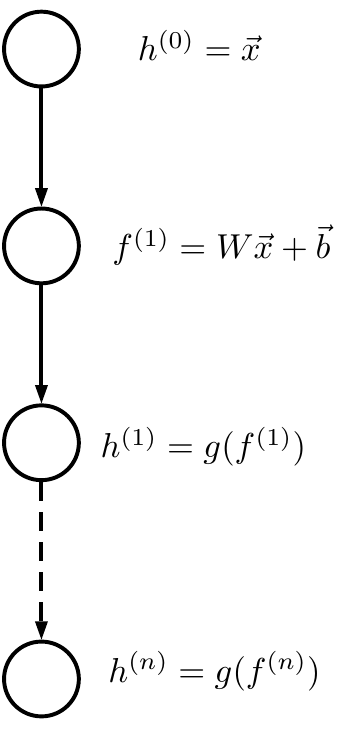}\hspace{3cm}
\includegraphics[scale=.7]{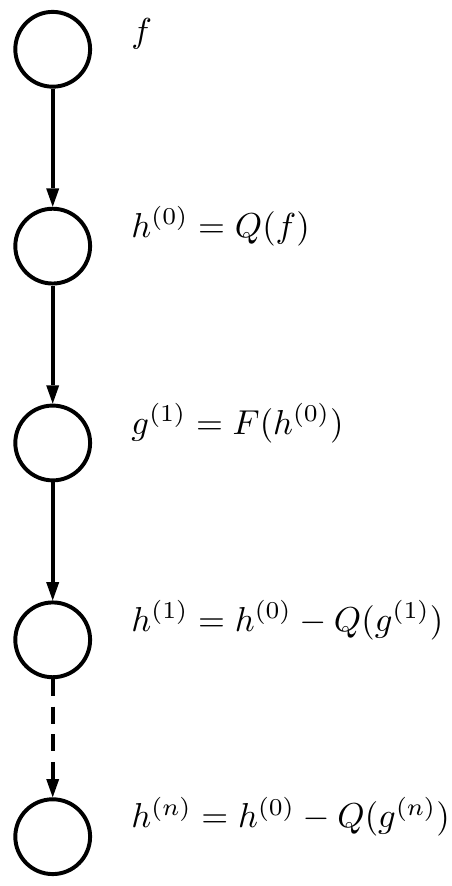}
\caption{Comparison of feedforward network and the iteration scheme for solving wave equations.}
\label{figcomp}
\end{figure}

By comparison, we can construct a simple neural network (MLP) from the iteration scheme for solving the wave equation: we replace the linear mapping by the solution operator of the linear wave equation and the activation function by $F(t, x, u).$ In case $F(t, x, u)$ is a polynomial of $u$ 
\beqq\label{eqpoly}
F(t, x, u) = \sum_{k = 2} ^N a_k(t, x) u^k = a_2(t, x) u^2  + \cdots + a_N(t, x)u^N,
\eeqq
we see that the parameters consists of the sound speed $c(x)$ (linear parameters) and the coefficients $a_k(t, x)$ of $F$ (nonlinear parameters). We shall denote 
\[
\theta = \{c(x), a_2(t, x), \cdots, a_N(t, x): (t, x)\in I(p_-, p_+)\}
\] 
In this setting, we see that all the layers and parameters in MLP have concrete meaning: the layers represents the propagation and nonlinear interactions in the solution and the parameters represent the significance of the nonlinearities. See Figure \ref{figcomp}.  

We state the network more precisely. Let $h^{(-1)} = f$ be the source term of the wave equation. This is the input data for the network. We set  
\[
h^{(0)} = Q(f)
\]
as the first layer which is just the linearized solution. Let $F(t, x, z)$ be a smooth function which is the activation function now. So we get 
\[
g^{(1)}  = F(t, x, Qh^{(0)}).
\]
It is better that we think of this as the hidden layer. Apply the linear operation to get 
\beqq\label{net1}
h^{(1)}  = h^{(0)} - Q F(t, x, Qh^{(0)}).
\eeqq
 This defines the iteration scheme and generates the second layer of the network. We then continue to obtain the network and obtain the output $h^{(N)}$. According to Proposition \ref{propite}, we immediately obtain 

\begin{theorem}\label{thm1}
Consider the inverse problem for wave equations with sources in Section \ref{sec-inv}. Assume 
\begin{enumerate}
\item $(f, u)\in \mcd_{sour}$ and $\|f\|_{H^s(V)}\leq \eps, s>1$ for $\eps$ sufficiently small.
\item $F$ is a polynomial function as   \eqref{eqpoly}.
\end{enumerate}   
Let $h^{(N)}$ be the approximation function generated from the network of depth $N$ using iteration \eqref{net1}. Then there exists parameter functions $\theta$ such that 
\[
\|u - h^{(N)}\|_{H^{s+1}(V)} \leq C_N \eps^{N+1}
\]
for some constant $C_N >0$ independent of $(f, u)$.
\end{theorem}

We remark that the theorem suggests that a deeper network produces better approximations and we have quantitative estimates to show this. To find the parameters which are directly related to $c$ and $F$, one solves the optimization problem with cost function $J = \|u - h^{(N)}\|^2_{H^{s+1}(V)}$ on a training set. 
A remarkable feature of deep neural network is the universal approximation property which allows one to determine the function $F$ without aa a-priori model (such as polynomials). This can be adapted to the network as follows. 

Again, we let $h^{(-1)} = f$ be the source term of the wave equation and set 
\[
h^{(0)} = Q(f)
\]
as the first layer. For $h^{(0)}$, we introduce $K$ units and apply affine transformations to get 
\beqq\label{eqappfun}
q^{(1)} = \sum_{k = 1}^K \gamma_k g( \alpha_k t + \beta_k, A_k x + B_k, a_k h^{(0)} + b_k)
\eeqq
where $g(t, x, u)$ is an activation function to be specified below and $\gamma_k, \alpha_k, \beta_k, a_k, b_k$ are constants and $A_k, B_k$ are constant matrix and vector. This step is supposed to approximate  $F(t, x, u)$. 
 Apply the linear operation to get 
\beqq\label{net2}
h^{(1)}  = h^{(0)} - Q(q^{(1)}).
\eeqq
This completes the first step and defines the iteration scheme. We then continue  to obtain the output $h^{(N)}$. The new parameter set is $\theta = \{c(x), a_k, b_k, \alpha_k, \beta_k,\gamma_k, A_k, B_k, k = 1, 2, \cdots, K\}$. Using Corollary \ref{thmuni1}, \ref{thmuni2}, we immediately obtain
 \begin{theorem}
Consider the inverse problem for wave equations with sources.  Assume 
\begin{enumerate}
\item $(f, u)\in \mcd_{sour}$ and $\|f\|_{H^s(V)}\leq \eps, s>1$ for $\eps$ sufficiently small.
\item  $F(t, x, u)$ is a smooth function with $F(t, x, 0) = F_u(t, x, 0) = 0.$
\end{enumerate}   
 Consider the network defined by iteration \eqref{net2} with activation function $g: \mbr^3\times \mbr \rightarrow \mbr$  in $C^l(\mbr^4)$ with non-negative integer $l \leq s$. Then there exists $K > 0$ and parameters $\theta$  such that 
\[
 \|u - h^{(N)}\|_{H^{l+1}(V)} \leq C_N \eps^{N+1},
 \]
 where $C_N$ is a constant independent of $f$ and $u.$ 
\end{theorem}

Again, one can solve the optimization problem on a training set to obtain the parameters, which further give approximations of $F(t, x, u)$ following \eqref{eqappfun}. 

It is important to realize that there are losses in the regularity of the estimates and this is essential for understanding our construction. For example, the ReLU activation function $\rho(x) = x_+$ is $C^0$. So the network only approximates the solution $u$ in $H^1$ norm. Obviously, $\rho(x)$ introduces new singularities to the network. Although $Q$ is a linear operator, it is non-local. Thus the new singularities might be propagated to other units. This issue does not show up in usual deep neural networks or the scattering network of Mallat \cite{Ma1}. In fact, the added singularity should help solving image classification problems from the singularity point of view, but not for our problem. 

Another issue is that Corollary \ref{thmuni1}, \ref{thmuni2} does not provide any estimate on the number of unit. In fact, to keep up with the $\eps^{N+1}$ error, the approximation error of $F(t, x, u)$  from the MLPs should be within $\eps^{N+1}$ instead of $\eps$. Thus one is not making good use of the nonlinearity. Roughly speaking, we think of the ``features'' in this inverse problem as the $H^s$ or $C^r$ singularities of the solution (or more precisely the wave fronts in phase space). This is similar to the ``edges'' in  images.   An important phenomena in nonlinear wave propagation is that nonlinear interactions of waves could produce new waves. This has been observed in physical applications and studied mathematically known as the nonlinear interaction of singularities and propagation of singularities for wave operators.  
 
Our next goal is to develop a neural network in Section \ref{sec-con} -- Section \ref{sec-app} which specifically addresses these issues.

 \section{Nonlinear interactions of conormal waves}\label{sec-con}
Conormal distributions have simple wave front sets and have been proven to be useful for analyzing wave interactions. In fact, otherwise the singularities generated from the nonlinear interactions could be rather complicated as shown by Beals example, see \cite{Bea}. 

We review Lagrangian distributions from H\"ormander \cite{Ho3, Ho4}. Let $X$ be a $n$ dimensional smooth manifold and $\La$ be a smooth conic Lagrangian submanifold of $T^*X\backslash 0$. We denote by $I^\mu(\La)$ the Lagrangian distribution of order $\mu$ associated with $\La$. In particular, for $U$ open in $X$, let $\phi(x, \xi): U\times \mbr^N \rightarrow \mbr$ be a smooth non-degenerate phase function that locally parametrizes $\La$ i.e. 
\beq
\{(x, d_x\phi): x\in U, d_\xi \phi = 0\} \subset \La.
\eeq
Then $u\in I^{\mu}(\La)$ can be locally written as a finite sum of oscillatory integrals
\beq
\int e^{i\phi(x, \xi)} a(x, \xi) d\xi, \ \ a\in S^{\mu + \fnf - \frac{N}{2}}(U\times \mbr^N),
\eeq
where $S^\bullet(\bullet)$ denotes the standard symbol class, see \cite[Section 18.1]{Ho3}.  For $u\in I^\mu(\La)$, we know that the wave front set $\WF(u)\subset \La$ and $u\in H^s(X)$ for any $s< -\mu-\fnf$. The principal symbol of $u$ is well-defined in $S^{\mu + \fnf}(\La; \Omega^\ha)/S^{\mu + \fnf- 1}(\La; \Omega^\ha)$, where $\Omega^\ha$ denotes the half-density bundle on $\La$. See Section 25.1 of \cite{Ho4}. For our problem, we can trivialize the bundle in local coordinates. 

For a submanifold $Y\subset M$, we denote $I^{\mu}(Y) = I^\mu(N^*Y)$, which are called conormal distributions to $Y$. In local coordinates $x = (x', x''), x'\in \mbr^k, x''\in \mbr^{n-k}$ such that $Y = \{x' = 0\}$. Let $\xi = (\xi', \xi'')$ be the dual variable, then $N^*Y = \{x' = 0, \xi'' = 0\}$. We can write $u\in I^\mu(Y)$ as
\beq
u = \int e^{ix'\xi'} a(x'', \xi')d\xi', \ \ a\in S^{\mu + \fnf - \frac{k}{2}}(\mbr^{n-k}_{x''}; \mbr^k_{\xi'}).
\eeq
In this case, the principal symbol is 
\beq
\sigma(u) = (2\pi)^{\fnf - \frac{k}{2}} a_0(x'', \xi')|dx''|^\ha |d\xi'|^\ha,
\eeq
where  $a_0 \in S^{\mu + \fnf - \frac{k}{2}}(\mbr^{n-k}_{x''}; \mbr^k_{\xi'})$ is such that $a - a_0 \in  S^{\mu + \fnf - \frac{k}{2}- 1}(\mbr^{n-k}_{x''}; \mbr^k_{\xi'})$. See \cite[Section 18.2]{Ho3}.  

Using four conormal waves and asymptotic analysis with multiple parameters, we can identify the leading terms in the solution that contains the new wave. This idea is introduced in Kurylev-Lassas-Uhlmann \cite{KLU1} and further developed in Lassas-Uhlmann-Wang \cite{LUW1}. We again consider the polynomial nonlinear function 
\[
F(t, x, u) = a(t, x) u^2 + b(t, x) u^3 + c(t, x) u^4
\]
We refine the iteration method in Section \ref{sec-ite} by introducing four small parameters to locate the nonlinear interactions. Let $f_i, i = 1, 2, 3, 4$ be compactly supported and set 
 \[
 f = \sum_{i = 1}^4 \eps_i  f_i
 \]
We let $v_i = Q(f_i), i = 1, 2, 3, 4$ be the linearized solution. Here, we shall assume that $v_i\in I^{\mu}(N^*Y_i)$ where $Y_i$ are codimension one submanifolds of $\mcm.$ These are called distorted plane waves, see \cite{KLU1, LUW1} for the details of construction in different context. With the source $f$, the linearized solution of $u$ is  
\[
v = \sum_{i = 1}^4  \eps_i v_i. 
\]
Let $u$ be the solution of the nonlinear equation and we use the  iteration scheme in Section \ref{sec-ite} to get 
\beq
\begin{gathered}
u  =  v + \sum \eps_i \eps_j Q(av_iv_j) + \mcr
\end{gathered}
\eeq
where the remainder term $\mcr = \sum_{i = 1}^4O(\eps_i^2)$ and the summation is over $i, j = 1, 2, 3, 4.$ In this approach, self-interactions of linearized waves are not considered. We iterate another two times to obtain  
\beq
\begin{gathered}
u  =  v - Q(F(v - Q(F (v - Q(F((v - Q(F(u))))))))\\
 =  v + \sum_{i, j}Q(av_i v_j) +  \sum_{i, j, k}\eps_i \eps_j \eps_k [Q(bv_i v_j v_k) + 2  Q(av_iQ(av_jv_k))] + \\
  \eps_1\eps_2\eps_3\eps_4\sum_{i, j, k, l} [Q(cv_iv_jv_kv_l) + Q(av_iQ(bv_jv_kv_l)) +   Q(bv_iv_jQ(av_k v_l)) + Q(av_iQ(av_jQ(av_kv_l)))] + \mcr. 
\end{gathered}
\eeq
We observe that the $\eps_i\eps_j$ terms reflects the interaction of two waves $v_i, v_j$, the $\eps_i\eps_j\eps_k$ terms the interaction of three waves and we are particularly interested in the $\eps_1\eps_2\eps_3\eps_4$ terms. It is worth noting that these terms can be obtained from 
\beqq\label{eqeps4}
\p_{ \eps_1}\p_{\eps_2}\p_{\eps_3}\p_{\eps_4} u|_{ \eps_1 = \eps_2 =  \eps_3 = \eps_4  = 0}. 
\eeqq

Suppose $Y_i, i = 1, 2, 3, 4$ intersect at a point $q\in I(p_-, p_+)$ transversally. The 
 of work \cite{KLU1, LUW1} shows that \eqref{eqeps4} at $q$ contains new singularities which are conormal to $T_q^*\mcm \backslash 0$. In other words, the term contains a point source. The  singularity can be propagated back to the region $V$ hence are observable in the data. Moreover, the leading order terms of the symbol of the conormal distributions are determined and they can be expressed in terms of the linear and nonlinear coefficients of the wave equation. The conclusion is that given all data $(f, u) \in \mcd_{sour}$, one can determine these coefficients in many cases up to diffeomorphisms, see \cite{KLU1, LUW1} for details. For illustration, we formulate a simple version of the uniqueness result below. 

\begin{theorem}\label{thmuni}
Let $c_1(x), c_2(x)$ be two smooth functions on $\mbr^3$ and let $g_1, g_2$ be associated Lorentzian metric. Let $V$ be a neighborhood of time like geodesics $\hat \mu_i \subset \mcm$. Let $-1<s_-<s_+ < 1$ and $p_i^\pm = \hat \mu_i(s_\pm)$. Consider the nonlinear wave equation $i = 1, 2$
\beq
\begin{gathered}
(\p_t^2 + c^2_i(x)\lap) u(t, x) + F_i(t, x, u(t, x)) = f(t, x), \quad (t, x) \in \mcm(T) \\
u(t, x) = 0, \quad (t, x)\in \mcm(0),
\end{gathered}
\eeq
where $F_i(t, x, u)$ are smooth such that $\p_u^k F(t, x, 0) \neq 0, x\in \mcm$ for some $k\geq 2.$ Assume that for $\delta$ sufficiently small, the data set 
\beq
\begin{split}
\mcd_{sour}^i  = & \{(f, u|_V): f\in C^4_0(V),\|f\|_{C^4} < \eps, \\
&\text{$u$ is the solution of nonlinear wave equation}\}, \quad i = 1, 2
\end{split}
\eeq
are the same. Then we have $c_1(x) = c_2(x)$ on $I(p_1^-, p_1^+) = I(p_2^-, p_2^+)$ and 
\[
\p_u^k F_1(t, x, 0) = \p_u^k F_2(t, x, 0), \quad k\geq 4. 
\]
\end{theorem}
\bpf
Because $c_i(x)$ does not depend on $t$, we consider the linearized problem and apply Tataru's unique continuation result to conclude that $c_1 = c_2$  on $I(p_1^-, p_1^+) = I(p_2^-, p_2^+)$. The determination of $F_i$ follows from \cite[Theorem 1.3]{LUW1}.  
\epf

We remark that if $c$ depends on $t$ or more generally one works with a globally hyperbolic Lorentzian metric $g$, Tataru's unique continuation result does not apply. One needs the full analysis in \cite{KLU, LUW1} and the determination is unique up to a conformal diffeomorphism in general. We also remark that the results are further applied to the Einstein equations coupled with scalar field equations or Maxwell equations, see \cite{KLU1, LUW2, UW}.

\section{The linear wave propagation}\label{sec-lin}
We consider linear variable coefficient wave equation 
\[
Pu = (\p_t^2 + c^2(x)\lap)u = f 
\] 
The fundamental solution $Q = P^{-1}$ is well understood. For our purpose, we need the microlocal structure of $Q$. In Melrose-Uhlmann \cite{MU},  a full symbolic construction was carried out and the Schwartz kernel $K_Q$ of $Q$ is found to be a paired Lagrangian distribution. We recall that for two Lagrangians $\La_0, \La_1 \subset T^*X$ which intersect cleanly at a codimension $k$ submanifold i.e. 
\beq
T_p\La_0\cap T_p\La_1  = T_p(\La_0\cap \La_1),\ \ \forall p\in \La_0\cap \La_1,
\eeq
the paired Lagrangian distribution associated with $(\La_0, \La_1)$ is denoted by $I^{p, l}(\La_0, \La_1)$.  For $u\in I^{p, l}(\La_0, \La_1)$, we know that $\WF(u)\subset \La_0\cup \La_1$. Microlocally away from the intersection $\La_0\cap \La_1$, $u \in I^{p+l}(\La_0\backslash \La_1)$ and $u \in I^p(\La_1\backslash \La_0)$ are Lagrangian distributions on the corresponding Lagrangians. 

Let  $\mcp(t, x, \tau, \xi) = |\tau|^2 - c^2(x)|\xi|^2, (t, x, \tau, \xi) \in T^*\mcm$  
 be the principal symbol of $P$.  Let $\Sigma$ be the characteristic set 
\beq
\Sigma  = \{(t, x,\tau, \xi)\in T^*\mcm: \mcp(t, x, \tau, \xi) = 0\}.
\eeq
The Hamilton vector field of $\mcp$ is denoted by $H_\mcp$ and in local coordinates
\beq
H_\mcp = \sum_{i =1}^4( \frac{\p \mcp}{\p \zeta_i}\frac{\p }{\p z_i} - \frac{\p \mcp}{\p z_i}\frac{\p }{\p \zeta_i}), \quad z = (t, x_1, x_2, x_3), \zeta = (\tau, \xi_1, \xi_2, \xi_3). 
\eeq
The integral curves of $H_\mcp$ in $\Sigma$ are called null-bicharacteristics. Let $\diag = \{(z, z')\in \mcm \times \mcm: z = z'\}$ be the diagonal and denote
\beq
N^*\diag = \{(z, \zeta, z', \zeta')\in T^*(\mcm \times \mcm)\backslash 0: z = z', \zeta' = -\zeta\}
\eeq 
  the conormal bundle of $\diag$ minus the zero section. We let $\La_{c}$ be the Lagrangian submanifold in $T^*(\mcm\times \mcm)$ obtained by flowing out $N^*\diag\cap \Sigma$ under $H_\mcp$. Here, we regard $\Sigma, H_\mcp$ as objects on product manifold $T^*\mcm \times T^*\mcm$ by lifting from the left factor. More explicitly, 
\[
\La_c = \{(z, \zeta, z', \zeta')\in T^*(\mcm\times \mcm): \text{ $(z, \zeta)$ lies on a bicharacteristics from $(z', -\zeta')$}\}
\]
The canonical relation is denoted by 
\[
\La_c' = \{(z, \zeta, z', \zeta')\in T^*(\mcm)\times T^*(\mcm):  (z, \zeta, z', -\zeta')\in \La_c\}
\]
We also call the map $S(z', \zeta') = (z, \zeta)$ if $(z, \zeta, z', \zeta')\in \La_c'$ the canonical relation. This map can be found explicitly by solving the Hamilton field equations. Let $\gamma(s) = (\alpha(s), \beta(s)): [0, \infty) \rightarrow \mcm\times \mbr^4$ be the null-bicharacteristics from $(z', \zeta')$. Then we have 
\beqq\label{eqode1}
\begin{gathered}
\frac{d\alpha(s)}{ds} = \frac{\p \mcp}{\p \zeta}, \quad \frac{d\beta(s)}{ds} = - \frac{\p \mcp}{\p z},\\
\alpha(0) = z', \quad \beta(0) = \zeta'.
\end{gathered}
\eeqq
Then $S(z', \zeta') = \gamma(s_0)$ where $\alpha(s_0) = z$. 
It is shown in \cite{MU} that for  linear differential operator $P$, the causal inverse $Q \in I^{-\frac{3}{2}, -\ha}(N^*\diag, \La_{c})$ is such that $PQ = \id$ on $\mce'(\mcm)$.   Also, from \cite[Prop.\ 5.6]{DUV}, we know that $Q: H_{\comp}^{s}(\mcm)\rightarrow H^{s+1}_{\loc}(\mcm)$ is continuous for $s\in \mbr$. 

We need the leading order singularities in $Q$. We follow the parametrix construction in \cite[Prop.\ 6.6]{MU}, see also \cite[Section 5.1]{Du}.  The conditions (6.1)-(6.6) of \cite{MU} are satisfied, thus the flow out of $\p \La_c$ under $H_\mcp$ is an embedded Lagrangian submanifold with boundary.  We look for $Q_0 \in I^{-\frac 3 2, -\ha}(N^*\diag, \La_c)$ to solve $PQ_0 - \id = 0$ with errors of  lower orders . First we have 
\[
\sigma(\id) = \mcp(z, \zeta) \sigma(Q_0)|_{N^*\diag}. 
\]
So on $N^*\diag$, we have $\sigma(Q_0) = \mcp(z, \zeta)^{-1}$. Then from \cite[Theorem 4.13]{MU}, we obtain (non-zero) initial condition of $\sigma(Q_0)$ on $\La_c\cap N^*\diag$. We solve on $\La_c$ 
\beqq\label{eqode2}
(i\mcl_{H_\mcp}+ \mcp_{sub} )\sigma(Q_0) = 0
\eeqq
where $\mcl$ denotes the Lie derivative acting on half density factors and $\mcp_{sub}$ is the subprincipal  
\[
\mcp_{sub} = -\frac{1}{2i}\sum \frac{\p^2 \mcp}{\p z_j \p \zeta_j}.
\]
Along null bicharacteristics from $(z', \zeta')$ to $(z, \zeta),$ the equation is a transport equation and we get the solution $\sigma(Q_0)(z, \zeta, z', \zeta')$ which is non-vanishing.  Using the canonical relation, we can write it as 
\[
\sigma(Q_0)(S(z', \zeta'); z', \zeta')
\]
So we find $Q_0\in I^{-\frac 32, -\ha}(N^*\diag, \La_c)$ such that 
\[
Q - Q_0 \in I^{-\frac{5}{2}, -\ha}(N^*\diag, \La_c). 
\]
Using the $L^2$ estimates of FIOs with paired Lagrangian kernel see \cite[Theorem 3.3]{GrU}, we obtain $Q-Q_0: H^s_{comp}(\mcm) \rightarrow H_{loc}^{s+2}(\mcm)$.

\section{Estimates of nonlinear effects}\label{sec-non}
Suppose $u\in H^{s+1}(\mcm)$ and $F(t, x, u)$ is smooth in $t, x, u$. Because in the linear wave propagation we only concerned the leading order singularities in $u$, we actually have $F(u) = F(w + R)$ where $u = w + R$ and $R\in H^{s+2}$.  In this section, we show that $F(u)$ can be approximated by a function $\tilde F(w)$  with difference difference $F(u) - \tilde F(w) \in H^{s+2}$ and such that $\tilde F(w)$ captures the nonlinear effects. Actually, we shall work in the phase space and make use of Bony's paraproducts. The approximation function $\tilde F$ is related to convolutional neural networks. 

We start with the dyadic decomposition of Coifman and Meyer \cite{CM}.  
 For $K > 1$ fixed, we set 
\[
\mcc_p = \{\xi\in \mbr^n: K^{-1} 2^{p} \leq |\xi| \leq K 2^{p+1}\},  
\]
 and $\mcc_0 = \{\xi\in \mbr^n:  |\xi| \leq 2K\}$.  Then $\{\mcc_p\}_{0}^\infty$ form an open covering of $\mbr^n.$  Let $\psi_j$ be a partition of unity 
\[
1 = \sum_{j = 0}^\infty \psi_j(\xi), \quad \psi_j \in C_0^\infty, \quad \supp \psi_j \subset \mcc_j.
\]
Actually, one can begin with $\psi_0(\xi)$ which is equal to $1$ for $|\xi|\leq K$ and 0 for $|\xi|> 2K$. Then set $\Psi_j(\xi) = \psi_0(2^{-j}\xi)$ and set $\psi_j(\xi) = \Psi_{j}(\xi) - \Psi_{j-1}(\xi)$. For any $u\in \mcs'(\mbr^n)$, the Paley-Littlewood decomposition of $u$ is
\beq
\{u_p\}_{0}^\infty \text{ where } u_p =\mcf^{-1} (\psi_p(\xi)\hat u(\xi)), p = 0, 1,2, \cdots.
\eeq
We have $u = \sum_{p = 0}^\infty u_p$ in the topology of $\mcs'(\mbr^n)$, see e.g.\ \cite{Tay}. 

We recall that Sobolev and H\"older functions can be characterized using Payley-Littlewood decompositions. The Sobolev space $H^{s}(\mbr^n)$ is defined as 
\[
H^{s}(\mbr^n) = \{u\in \mcs'(\mbr^n): (1+ |\xi|^2)^{s/2}\hat u(\xi) \in L^2(\mbr^n)\}
\]
with norm
\[
\|u\|_{H^s} = \|(1 + |\xi|^2)^{s/2}\hat u(\xi)\|_{L^2}. 
\]
Then $u\in H^s(\mbr^n)$ if and only if $u = \sum_{p = 0}^\infty u_p$ where $\hat u_p$ are supported in $\mcc_p$ and  
\[
\|u_p\|_{L^2}\leq c_p 2^{-ps}, \quad \{c_p\} \in l^2.
\]
Consider the H\"older space $C^\alpha(\mbr^n), \alpha > 0$ non-integer, equipped with the norm 
\[
\|u\|_{C^\alpha} = \sum_{|\la|\leq [\alpha]} \|\p^\la u\|_{C^\beta}.
\]
When $\alpha$ are integers, it is necessary to use the Zygmund space $C_*^{\alpha}(\mbr^n)$.  
In particular, $C_*^\alpha = C^\alpha$ if $\alpha$ is not integer. Otherwise, $C^\alpha \subset C^\alpha_*$. The characterization is that $u\in C^\alpha_*(\mbr^n)$ if and only if $u = \sum_{p = 0}^\infty u_p$ where $\hat u_p$ are supported in $\mcc_p$ and  
\[
\|u_p\|_{L^2}\leq c 2^{-p\alpha}. 
\] 

Let $a\in C^r(\mbr^n)$ and $f\in H^{s}(\mbr^n)$.  
The paraproducts of $a$ and $f$, introduced by Bony \cite{Bo}, is 
\beqq\label{defpara}
T_a f = \sum_{k\geq 1} (\Psi_{k-1}(D)a) (\psi_{k+1}(D)f) = \sum_{p = 2}^\infty \sum_{q =0}^{p-2} a_q f_p.
\eeqq
where $\Psi_k(\xi) = \sum_{j = 0}^k \psi_j(\xi).$ If we denote $\mcb_p = \{\xi\in \mbr^n:  |\xi| \leq K 2^{p+1}\}$, then $\Psi_k$ is supported in $\mcb_k.$ 
Using the characterization of $H^s, C^r$ functions, we see that $T_af\in H^{s}(\mbr^n)$ and 
\[
af = T_a f+ R, \quad R\in H^{s+ r}(\mbr^n).
\]
So the difference is a more regular term for $r > 0$. Furthermore, we have that if $u\in C^r\cap H^s, r, s > 0$ and $F(u)$ is smooth in $u$, then 
\[
F(u) = T_{F'(u)} u + R, \quad R\in H^{s+ r}.
\]
See \cite[Proposition 3.2.C]{Tay}. We remark that the paraproduct does not throw away all nonlinear effects, which is evident from the definition \eqref{defpara}.  

There are several equivalent variants of paraproducts, see \cite{Tay}. The one convenient for our purpose is to introduce a convolution kernel in the phase space which is also done in Bony \cite{Bo}. Choose $\chi \in C^\infty(\mbr^n\times \mbr^n)$, homogeneous of degree 0 outside a compact set such that $\chi(\xi, \eta) = 0$ for $|\xi| > \ha |\eta|$ and $\chi(\xi, \eta) = 1$ for $|\xi| < |\eta|/16$ and $|\eta| > 2$. Then the paraproduct can be written as 
\beq
\begin{gathered}
T^\chi_a f(x) = (2\pi)^{-n}\int e^{ix\xi} \chi(\xi - \eta, \eta) \hat a(\xi - \eta)\hat f(\eta)d\eta d\xi\\
 = (2\pi)^{-n} \int e^{i x (\xi+ \eta)} \chi(\xi, \eta) \hat a(\xi)  \hat f(\eta) d\xi d\eta. 
\end{gathered}
\eeq
We see that $\widehat{T^\chi_a f}$ is a convolution of $\hat a, \hat f$ with kernel $\chi$. 
Here, we emphasized the dependence on $\chi.$ We also use the notation $T_a^\chi f = T^\chi(a; f)$. \\
 
We use paraproducts to construct a network for approximating composite functions $F(u), u\in H^s(\mbr^n)$. (Here, $F$ is only a function of $u$ not $x$.) Let $h^{(0)} =  u$ be the first level of the network. Then we perform affine transformations and use paraproducts as the activation function to get 
\[
h^{(1)}(u) = T^\chi(u; a_1u + b_1)
\]
where $a_1, b_1$ are constants. Then we continue to get the $n$-th step   
\beqq\label{netconv}
h^{(n)}(u) = T^\chi(u; a_n h^{(n-1)} + b_n)
\eeqq
We remark that by taking the Fourier transform,  the network is a convolutional network with kernel $\chi$ and $\hat h^{(n)}$ and the network  provides an approximation of $\widehat F(u)$ in the phase space. These two point of views will be used interchangeably below. 

We analyze the difference $F(u) - h^{(n)}(u)$. We first prove that the error terms consist of a spacial error which is controlled by the nonlinearity of $F$ and a phase space error term which is controlled by the regularity of $u$. 
For $R>0$, let $\Psi_R(\xi), \xi\in \mbr^n$ be a smooth cut-off function such that $\Psi_R(\xi) = 0$ if $|\xi| < R$ and $\Psi_R(\xi) = 1$ if $|\xi|> 2R$. We denote $\Psi_R(D)$ be the pseudo-differential operator with symbol $\Psi_R$. 

\begin{prop}\label{propapp}
We assume that 
\begin{enumerate}
\item $u\in H^{s}(\mbr^n)\cap C^r(\mbr^n), r, s > 0$ with $\|u\|_\infty <\eps$;
\item $F(u)$ is a smooth function of $u$.
\end{enumerate}
Then there exists constant parameters $a_i, b_i, i = 1, \cdots, N$ such that for the $h^{(N)}(u)$ obtained in \eqref{netconv}, we have $F(u) - h^{(N)}(u) = R_{sp} + R_{ph}$ where $R_{sp}\in H^{s}, R_{ph}\in H^{s+r}$ and 
\[
\|R_{sp}\|_{H^s} < C_F\eps^{N+1}, \quad \|\Psi_R(D) R_{ph}\|_{H^{s}} = O(R^{-r}). 
\] 
Here $C_F$ is a constant such that $\sup_{\|u\|_\infty < \eps}|\p_u^{N+1}F(u)| \leq C_F$. 
\end{prop}

\bpf
The proof is straightforward. First we use Taylor expansion of $F$ based at $0$ 
\[
 F(u) = \sum_{n = 0}^N \frac{\p_u^n F(0)}{n!} u^n + R_{sp}, \quad |R_{sp}| \leq C_F |u|^{N+1}. 
\]
Next, we let $p(u) =  \sum_{n = 0}^N \frac{\p_u^n F(0)}{n!} u^n$ and rewrite it as 
\[
p(u) = a_nu(\cdots a_3u(a_2u(a_1 u + b_1) + b_2)) + b_3\cdots) + b_n 
\]
where $a_i, b_i$ are constants related to $\frac{\p^n F(0)}{n!}$. Now we replace the products by paraproducts. First of all, 
\[
a_2u(a_1 u + b_1) + b_2 = a_2T^\chi(u; a_1 u + b_1) + b_2 + R_2
\]
where $R_2 \in H^{s+r}$. Next, we get
\[
a_3(a_2u(a_1 u + b_1) + b_2) + b_3 = a_3T^\chi(u; a_2T^\chi(u; a_1 u + b_1) + b_2) + a_3 T^\chi(u; R_2) + R_3
\]
where $R_3\in H^{s+r}$ and we also have  $ a_3 T^\chi(u; R_2) \in H^{s+2r}$.  Continuing this procedure, we get the function $h^{(N)}$ such that $p(u) - h^{(N)}(u) = R_{ph} \in H^{s+r}$. 
This finishes the proof. 
\epf

Next, we show the stability of the network with respect to regular perturbations. 
\begin{cor}\label{propreg}
We assume that 
\begin{enumerate}
\item $u, w\in H^{s}(\mbr^n)\cap C^r(\mbr^n), r, s > 0$ with $\|u\|_\infty,\|w\|_\infty <\eps$;
\item $u = w + R$ with $R\in H^{s+m}(\mbr^n), m>0$; 
\item $F(u)$ is a smooth function.
\end{enumerate}
Then there exists constant parameters $a_i, b_i, i = 1, \cdots, N$ such that for the $h^{(N)}(w)$ obtained in \eqref{netconv}, we have $F(u) - h^{(N)}(w) = R_{sp} + R_{ph}$ where $R_{sp}\in H^{s}, R_{ph}\in H^{s+r}$ and 
\[
\|R_{sp}\|_{H^s} < C_F\eps^{N+1}, \quad \|\Psi_R(D) R_{ph}\|_{H^{s}} = O(R^{-t}) 
\] 
where $t = \min(m, r)$ and $C_F$ is the same as in Prop.\ \ref{propapp}.
\end{cor}
\bpf
The spacial error is the same as in the previous proposition. So we consider 
\[
p(u) = \sum_{n = 0}^N \frac{\p_u^n F(0)}{n!} (w + R)^n =  \sum_{n = 0}^N \frac{\p_u^n F(0)}{n!}  w^n + R_0\]
where $R_0\in H^{s+m}$ because it is a finite sum of products of $w \in C^r$ and $R\in H^{m}$. This finishes the proof.
\epf

We remark that if $F(u)$ is such that $F(0) = F'(0) = 0$, we see from the proof that the estimate of $R_{ph}$ is actually $\|\Psi_R(D) R_{ph}\|_{H^{s}} = O(\eps^2 R^{-t})$.  Finally, we make a remark  in relation to the usual convolutional networks, see e.g. \cite{LeC, Ma2}. Let $u\in H^s\cap C^r$ and we consider the paraproduct of $u$: 
 \beq
\begin{gathered}
\widehat T^\chi_u u(x) = \int  \chi(\xi - \eta, \eta) \hat u(\xi - \eta)\hat u(\eta)d\eta 
\end{gathered}
\eeq
Now we use the Paley-Littlewood decomposition and consider for $M$ large
\[
u\simeq \sum_{p = 0}^M u_p  = \sum_{p = 0}^M \mcf^{-1}( \psi_p(\xi)\hat u(\xi)).
\] 
Then formally, we obtain
\beq
\begin{gathered}
\widehat T^\chi_u u(x) \simeq  \int \sum_{p = 0}^M \chi(\xi - \eta, \eta)\psi_p(\xi-\eta)\hat u(\xi-\eta)  \hat u(\eta)d\eta 
\end{gathered}
\eeq
One can think of $\chi^p(\xi, \eta) = \chi(\xi, \eta)\psi_p(\xi)$ as the convolutional kernel at different scales. 
Here, the other $\hat u(\xi-\eta)$  plays an important role because it captures the nonlinear effects. However, if one approximates this   $\hat u(\xi-\eta)$  on the support of $\chi^p$ using a set of parameters, one would obtain the usual convolutional network. Moreover, keeping the leading order terms at every step is in the same spirit as the max pooling operation (see \cite{LeC, Go}). These considerations, in some sense, show that if we interpret the above network as a conventional convolutional neural network, some of the parameters are related to the data themselves and some of them to the nonlinear functions. 

\section{Convolutional neural network in the phase space}\label{sec-conv}
We construct the network for solving the inverse problem in phase space. We return to the setup in Section \ref{sec-inv}. Let $V$ be an open relatively compact set of $\mcm$ and $f$ be the source function supported in $V$. This is the input data for the network. 

We first choose a finite open covering $U_i, i = 1, \cdots K$ of $I(p_-, p_+)$ where $p_\pm \in V$. 
Let $\diam(U)$ be the diameter of set $U\subset \mcm$ and we assume that $\diam(U_i) < \delta, i = 1, 2, \cdots, K$. We see that $K$ is at least $O(\delta^{-2}).$ We let $\phi_i$ be a partition of unity subordinated to $U_i $ 
\[
\phi_i \in C_0^\infty(U_i), \quad \sum_{i = 1}^K \phi_i = 1.
\]
For any $f\in H^s(\mcm)$, we write $f_i = \phi_i f  \in H^s(\mcm)$ and get $f = \sum_{i = 1}^K f_i. $ Then we define the $0$-th layer (input) of the network to be 
\[
h^{(0)} = \{\hat f_i\}_{i = 1}^K.
\]
In particular, the number of units for this level is $K$. When $U_i$ is taken sufficiently small, $\hat f_i$ is a good approximation of the wave front set of $f$ at $U_i$. 

\begin{figure}[htbp]
\centering
\includegraphics[scale = 0.7]{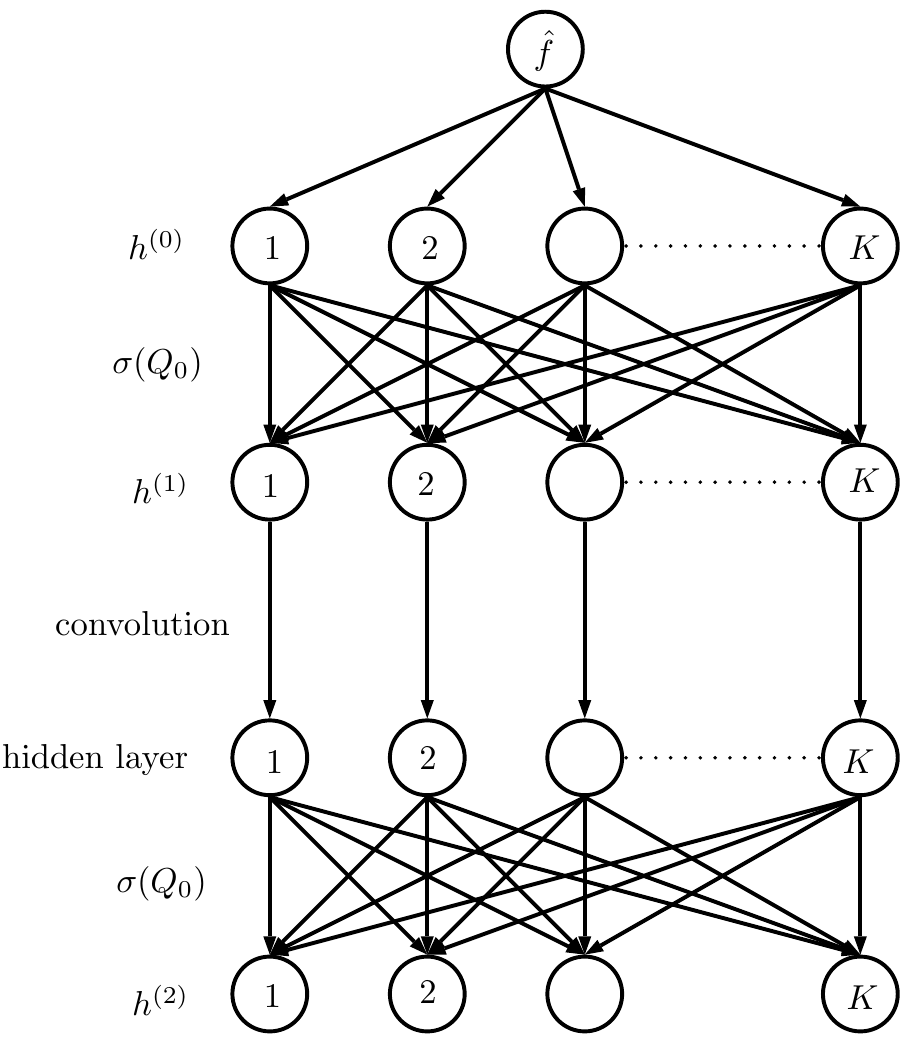}\quad \quad
\caption{Illustration of the network. $h^{(0)}$ is the input layer which we regard as the $0$-th layer. $h^{(1)}, h^{(2)}$ are the first and second layer of the network.  }
\label{figdnn}
\end{figure}

Next, we solve the wave equation from each $U_i$ to $U_j$. This means that we solve
\[
Pu = f_i \text{ in } \mcm
\]
and get $u|_{U_j}$. This makes sense if $U_i\cap J_+(U_j) \neq \emptyset.$ We remark that one can regard $u$ as the wave-packet generated by a point source if $U_i$ is sufficiently small. We shall use the leading term of $Q$ on $\La_c$ with principal symbol $\sigma(Q_0).$ For each pair $U_i, U_j$, we further decompose the operation as follows. Let $S_{ij}$ be an invertible matrix which is an approximation of the  canonical relation $S$. Then we set 
\[
h^{(1)}_{ij}(\zeta) = c_{ij}|\zeta|^{-1} h^{(0)}_i(S_{ij}\zeta), \quad \zeta\in \mbr^4, \quad i\neq j
\]
which we think of as an approximation of $\hat{\phi_iu}$ (to be justified later). If $U_{i}\cap J_+(U_j) = \emptyset$, we should take $c_{ij} = 0$. This step solves wave propagation. On $U_j$ itself, we solve using $Q_0$ on $N^*\diag$ so 
\[
h^{(1)}_{jj}(\zeta) = \frac{\gamma(\zeta)h^{(0)}_j(\zeta)}{|\tau|^2 - c^2_{jj}|\xi|^2}, \quad \zeta = (\tau, \xi) \in \mbr^4. 
\]
One can think of $c_{jj}$ as the constant wave speed on $U_j$ and $\gamma(\zeta)$ is a cut-off function away from the light-like directions. In particular, let $\mcp_{j}(\zeta) = |\tau|^2 - c^2_{jj}|\xi|^2|$. Then we take $\gamma(\zeta) =0$ when $\mcp_j(\zeta)< \delta$ and $\gamma(\zeta) =1$ if $\mcp_j(\zeta)> 2\delta$.  
 We collects the effects on each $U_i$ and let 
\[
 h^{(1)} = \{h^{(1)}_i\}_{i = 1}^K, \quad h^{(1)}_i = \sum_{j = 1}^K   h^{(1)}_{ij}.
\]
We shall take this as the first layer of the network. See Figure \ref{figdnn}.

To obtain the next layer, we need to take into account the nonlinear effects and find approximations of $\hat F$. Now we use the convolutional network constructed in Section \ref{sec-non} on each $U_i$ with parameter set $\theta_{ik} \doteq \{a_{ik}, b_{ik}\}, i = 1, 2, \cdots, K; k = 1, 2, \cdots, M$. We denote the obtained approximation function on $U_i$ by $v_i^{(1)}, i = 1, 2, \cdots, K$. Next, we take $v_{j}^{(1)}$ as the source to solve the wave equation in phase space as before to get $h^{(2)}_{ij}$ on $U_i$. Again, the terms are no longer supported on $U_i$ so we collect the terms on each $U_i$ to obtain the second layer
\[
h^{(2)}_i =  h^{(1)}_i - \sum_{j = 1}^M  h^{(2)}_{ij}
\] 
This layer collects the linear and quadratic effects in the solution. See Figure \ref{figdnn} for the illustration of the structure. 
Continue the procedure, we obtain the approximation function $h^{(M)}$ which only depends on the parameter set 
\[
 \Theta = \{c_{ij}, S_{ij}: i, j = 1, \cdots, K\} \cup \{\theta_{ik}: i =1,\cdots, K, k = 1, \cdots, M\}
 \]
  where we take $S_{ii}  = \id$ the identity. We denote the approximation by $h^{(M)}(f; \Theta)$.

To determine the parameters, we need to solve the optimization problem on training data with the cost function 
 \[
 J_i = \|\hat u_i - h^{(M)}_i(f; \Theta)\|^2, \quad i\in \mci,
 \]
 where the index set $\mci$ is such that $U_i \subset V, i\in \mci$ and $u_i = \phi_i u$ is supported on $U_i$. We shall specify the proper norm after the analysis in next section. We remark that  one can use any of these cost functions or combinations of them. 

The next section is devoted to the approximation properties of this network.

\section{The approximation theorem}\label{sec-app} 
\begin{theorem}\label{thmmain}
Consider the inverse problem for nonlinear wave equations with sources formulated in Section \ref{sec-inv}. Assume that
\begin{enumerate}
\item  $c(x), F(t, x, u)$ are smooth functions. $F(t, x, 0) = F_u(t, x, 0) = 0$. 
\item $(f, u) \in \mcd_{sour}$ and $\|f\|_{H^s(V)}\leq \eps < \eps_0, s>1$, where $\eps_0$ is as in Prop.\ \ref{propite}.
\end{enumerate}

Consider the convolutional network constructed in Section \ref{sec-conv} with depth $M\geq 0$ and $K\geq 1$ units for each level. Let $\mci$ be the index set so that $U_i\subset V$ and let $u_i = \phi_i u$. Assume that $\text{diam}(U_i)<\delta$ for some $\delta > 0$. 
Then there exist parameter sets $\Theta$ and $M, K$ such that the function $h^{(M)}(f; \Theta)$ generated by the neural network satisfy for $i\in \mci$
\beq
\begin{gathered}
\|\Psi_R(\zeta)\langle\zeta\rangle^{(s+1)/2}(\hat u_i(\zeta) - h_i^{(M)}(f; \Theta))\|_{L^2(\mbr^4)}
 \leq C(1+ \delta ) \eps^{M} 
 \end{gathered}
\eeq
where $R> R_0$ and the constant $C, R_0$ and $K$ depends on $M, \delta, \eps_0, c(x)$ and $F(t, x, u)$. 
\end{theorem}

We make few remarks before giving the proof. This theorem indicates that it is better to solve the optimization problem in the phase space and consider high frequency information. We shall see in the proof that the error comes from two sources. One is 
$\|u - u^{(n)}\|_{H^{s+1}} \leq C_n\eps^n$ from Prop.\ \ref{propite} where $C_n$ is found in the proof. The other one is 
\[
\|\Psi_R(\zeta)\langle\zeta\rangle^{(s+1)/2}(\hat u^{(n)}_i(\zeta) - h_i^{(M)}(f; \Theta))\|_{L^2(\mbr^4)}
 \leq C\delta  \eps^{n} 
\]
and we shall see that $C$ depends on $\eps_0, M$ and 
\[
\sup_{x\in I(p_-, p_+)}|c(x)|,  \quad \sup_{x\in I(p_-, p_+)}|\p_x c(x)|, \quad \sup_{(t, x)\in I(p_-, p_+), |s|< 1 }|\p_s^k \p_{(t, x)} F(t, x, s)|, \quad  k\leq M. 
\]
Here, $x\in I(p_-, p_+)$ means $x$ in the projection of $I(p_-, p_+)$ to $\mbr^3.$ By taking $\delta$ small (necessarily increasing the number of units $K$), we obtain better approximation results.  Indeed, when $\delta \rightarrow 0$, the set $U_i$ approaches to a point. Essentially what matters in the units of the network is just the wave front sets of $u$ so the estimates become more accurate. Finally, we remark that in view of the uniqueness result Theorem \ref{thmuni} and its proof, one can take the training data consisting of sufficiently many conormal waves that are supported on each $U_i, i\in \mci.$

\bpf[Proof of Theorem \ref{thmmain}]
Because $s > 1$, we know from Prop.\ \ref{propite} and Sobolev embedding that the solution $u\in H^{s+1} \subset C^r, r < s-1$. 
We start with the first level, that is $h^{(1)}$. This involves solving the wave equation using $Q_0.$ Recall the open covering $U_i, i = 1, 2, \cdots, K$ for $I(p_-, p_+)$ and $f_i = \phi_i f \in H^s$ are compactly supported on $U_i$. 

We first solve $P v = f_i$ away from $U_i$. From Section \ref{sec-lin}, we know that $v - Q_0(f_{i})  \in H^{s+2}$. Away from $U_i$, it suffices to consider $Q_0\in I^{-\frac 32}(\La_c\backslash N^*\diag)$. So we can write
\[
Q_0(f_{i})(z) = \int e^{i \phi(z, z', \theta)} a(z, z', \theta) f_{i}(z')dz' d\theta
\]
where $\phi$ is a homogeneous non-degenerate phase function that parametrizes the Lagrangian $\La_c$ locally near $(z, z')$, namely
\[
\La_c = \{(z, \zeta, z', \zeta') \in T^*(\mcm\times \mcm)\backslash 0: \zeta = \phi_z, \zeta' = -\phi_{z'}, \phi_{\theta} = 0\}
\]
and $a$ is a smooth function homogeneous of degree $-1$ in $\theta$ for $|\theta| > 1$. Consider $z\in U_j, z' \in U_i, i\neq j$ and choose constant $\tilde c_{ij}$ such that 
\[
|\tilde c_{ij}\langle \theta\rangle^{-1} - a(z, z', \theta)| \leq C \delta \langle\theta\rangle^{-1}
\]
Here, $C$ depends on the symbol $a$. Because the symbol is obtained by solving the transport equation  \eqref{eqode2} involving $c(x)$ and its first derivatives along null-bicharactersitics, from the stability of ODEs, we see that $C$ depends on $\|c\|_{C^1} $ on $I(p_-, p_+)$.  Then we have 
\beq
\begin{gathered}
\phi_i(z)Q_0(f_{i})(z)  - \phi_i(z)\int e^{i \phi(z, z', \theta)} \tilde c_{ij}\langle \theta\rangle^{-1} f_{i}(z')dz' d\theta = \mcr_1, \\
\text{where } \|\mcr_1\|_{H^{s+1}(U_i)} \leq C \delta \|f\|_{H^s(V)}\leq C\delta \eps. 
\end{gathered}
\eeq
Now we denote  
\beq
I(z) = \phi_i(z)\int e^{i \phi(z, z', \theta)} \tilde c_{ij}\langle \theta\rangle^{-1} f_{i}(z')dz' d\theta
\eeq
 and take the Fourier transform to get 
\begin{eqnarray*}
\hat I(\zeta) &\doteq& \int e^{-iz\zeta} e^{i \phi(z, z', \theta)} \phi_i(z) \tilde c_{ij}\langle \theta\rangle^{-1} f_{j}(z')dz' d\theta dz\\
& =& (2\pi)^{-4} \int  e^{-iz\zeta} e^{i \phi(z, z',  \theta)} e^{i z' \eta} \phi_i(z) \tilde c_{ij}\langle \theta\rangle^{-1} \hat f_{j}(\eta)dz' d\theta dz  d\eta.
\end{eqnarray*}
For this oscillatory integral, the phase function is 
\[
\Phi(z, z', \zeta, \theta, \eta) = -z\zeta + \phi(z, z', \theta) + z'\eta
\]
which is non-degenerate and homogeneous of degree one in $\zeta, \eta, \theta$. The critical points are 
\beq
\begin{gathered}
\Phi_{\theta} = \phi_{\theta} = 0, \quad \Phi_{\eta} = z' = 0, \quad \Phi_z = -\zeta + \phi_z = 0, \quad \Phi_{z'} = \eta + \phi_{z'} = 0. 
\end{gathered}
\eeq
Suppose $(z, \zeta; z', \zeta')\in \La_c'$ and we choose local coordinate so that $z' = 0$. 
Using stationary phase arguments (e.g.\ \cite[Prop.\ 1.2.4]{Du}), we obtain that 
\[
\hat I(\zeta)= c_{ij} \langle\zeta\rangle^{-1}\hat f_{j}(\zeta') + \mcr_1'
\] 
 with new parameters $c_{ij}$ and where $\langle \zeta\rangle^{(s+2)/2}\mcr_1'\in L^2$ and $\|\mcf^{-1} \mcr_1'\|_{H^{s+2}} = O(\eps)$.

Recall that $S(z, \zeta) = (z', \zeta')$. Let $S_{ij}$ be a $4\times 4$ matrix such that  $|S_{ij} \zeta - S(z, \zeta)| \leq C\delta |\zeta|$ for $z \in U_i, z'\in U_j, |\zeta| > 1$. Here, $C$ depends on $S$. But we know from Section \ref{sec-lin} that $S$ is the solution of ODEs \eqref{eqode1} with coefficients depending on $c(x)$ and its first derivatives. By the stability of ODEs, we see that $C$ depends on $\|c\|_{C^1}$ on $I(p_-, p_+).$ Then we get  
\[
\hat I(\zeta)= c_{ij}\langle \zeta\rangle^{-1} \hat f_{j}(S_{ij}\zeta) + \mcr_2 = h_{ij}^{(1)}(f; \Theta) + \mcr_2, \quad i\neq j. 
\]
Here, $\Theta$ is the collection of parameters of the network including $c_{ij}, S_{ij}$.  To estimate $\mcr_2$, we recall that $f_j \in H^{s}$ and we have for $k\leq s$ that
\beq
\begin{gathered}
 \langle \zeta \rangle^k |\hat f_j(\zeta + \delta\zeta) - \hat f_{j}(\zeta)| = \langle \zeta \rangle^k |\int  (e^{iz(\zeta + \delta \zeta)} - e^{iz \zeta }) f_j(z) dz| \\
\leq  C \delta \sum_{|\alpha| = k}  | \int \p_z^\alpha (e^{iz(\zeta + \delta \zeta)} - e^{iz \zeta }) f_j(z)dz |\leq C \delta \|f\|_{H^{k}}
\end{gathered}
\eeq
because $f_j$ is supported in $U_i \subset I(p_-, p_+)$. So the constant $C$ depends on the size of $I(p_-, p_+)$. Let $v_{i} = \phi_i v$.  
  Therefore, we proved that 
\beqq\label{eqres11}
\|\Psi_R(\zeta)\langle\zeta\rangle^{(s+1)/2}(\hat v_i - h_{ij}^{(1)}(f; \Theta))\|_{L^2} \leq C\delta \|f\|_{H^s(V)} + O(\eps R^{-1}) \leq C\eps(\delta  + R^{-1}) \leq C\eps \delta
\eeqq
if $R > 1/\delta$ is large enough. 

Next, we consider solving $Pv = f_j$ on $U_j$. We want to use $Q_0$ on $N^*\diag$ which is a pseudo-differential operator so we ignore the part on $\La_c$. So we introduce a microlocal cut-off $\Phi$ supported sufficiently close to $\La_c\cap N^*\diag$. In particular, we let $\chi(z, \zeta)$ be smooth in $T^*\mcm$ and $\chi(z, \zeta) = 1$ in $\mcp(z, \zeta)< \delta$ and $\chi(z, \zeta) = 0$ in $\mcp(z, \zeta) > 2\delta$. Then let $\tilde\chi(t)$ be a smooth cut-off function so that $\tilde\chi(t) = 1$ for $|t|< \delta$ and $\tilde \chi(t) = 0$ for $|t| > 2\delta$. Then we set $\Phi(z, \zeta, z', \zeta') = \chi(z, \zeta)\tilde\chi(|z-z'|+ |\zeta-\zeta'|)$. Because $\diam(U_i) < \delta$, we still have a $\delta$ order error. More precisely,
\beq
\phi_i(z) Q_0(f_j) - \phi_i(z) \int e^{i(z-z')\zeta}(1-\chi(z,\zeta))\frac{ f_j(z') }{\mcp(z,\zeta)}dz' d\zeta = \mcr_3
\eeq
where $\|\mcr_3\|_{H^{s+1}}\leq C\delta \|f\|_{H^{s}(V)}\leq C\delta \eps$ and $C$ depends on the symbol only. Then we estimate  
\beq
\phi_i(z) \int e^{i(z-z')\zeta}(1-\chi(z,\zeta) )\frac{ f_j(z') }{\mcp(z,\zeta)}dz' d\zeta - \phi_i(z) \int e^{i(z-z')\zeta} (1 - \chi(z, \zeta)) \frac{f_j(z')}{|\tau|^2 - c^2_{jj} |\xi|^2} dz' d\zeta = \mcr_4
\eeq
Using the same argument as we used for \eqref{eqres11}, we see that $\|\mcr_4\|_{H^{s+2}} \leq C\delta \|f\|_{H^s}\leq C\delta \eps$ if $|c^2_{jj} - c^2(x)| \leq C\delta$ on $U_i.$ So we proved \eqref{eqres11} for $v_j = \phi_j v.$

Now we consider the second layer $h^{(2)}$ and we need the nonlinear function $F(t, x, u)$. On each $U_i$, we write $F(t, x, u)$ in Taylor expansions
\[
F(t, x, u) = \sum_{j = 2}^M a_{j}^{(i)}(t, x) u^j + O(|u|^{M+1}). 
\]
Let $a_{ij}$ be constants so that $|a_{ij} - a^{(i)}_j(t, x)| < C\delta$ where $C$ is some  constant depending on the sup norm $|\p_{(t, x)} a^{(i)}_j(t, x)|_\infty, i = 1, 2, \cdots, K$. Thus for $p(u) = \sum_{j = 2}^M a_{ij} u^j$, we obtain that $|F(t, x, u) - p(u)| \leq C\delta |u|^2$.  Let $v_j^{(1)}$ be obtained  from $h^{(1)}_j$ in the network and $v_j$ be the solution of linearized wave equation. Using \eqref{eqres11}, we apply Prop.\ \ref{propapp}, Corollary \ref{propreg} and the remark after them to get that $\hat{p(v_j)} - v^{(1)}_j = \hat R_{sp} + \hat R_{ph}$ and they satisfy  
\beq
\|R_{sp}\|_{H^{s+1}} < C \delta \eps^{2}, \quad \| \Psi_R(D) R_{ph}\|_{H^{s+1}} \leq C \eps^2 \delta  R^{-r} 
\eeq
where $r < s-1$.  Next, we can apply the argument above to solve the linear wave equation $P w = v_{j}^{(1)}$ on each $U_i$ to get $w_i.$ Then we obtain the term $h^{(2)}_{ij}$ in the network and 
\[
\hat w_i - \sum_{j = 1}^K h^{(2)}_{ij} =  \mcr_5
\]
where 
\[
\|\Psi_R(\zeta)(1 + |\zeta|)^{(s+1)/2}\mcr_5\|_{L^2} \leq C\delta  \eps^2(1 + R^{-r})
\]
for $R$ large enough. Together with Proposition \ref{propite}, we complete the analysis for the iteration step in the network and obtain
\beq
\|\Psi_R(\zeta)\langle\zeta\rangle^{(s+1)/2}(\hat u_i(\zeta) - h_i^{(2)}(f; \Theta))\|_{L^2(\mbr^4)}
 \leq  C_2 \eps^2 + C \delta \eps^{2}(1 + R^{-r}) 
\eeq 
where $C_2$ is the constant in Prop.\ \ref{propite} which depends on $c$ and $F.$ 
The proof is finished by induction. 
\epf

Finally, we discuss the reconstruction of $c(x)$ and $F(t, x, u)$ on each $U_i$ from the parameters. For fixed $i = 1, 2, \cdots, K,$ consider the collection of $\theta_{ik}, k = 1, 2, \cdots M$ which are the parameter sets on $U_i$. From the construction of the network and the proof, it is easy to see that 
\[
p_i(u) = a_{iM}u(\cdots a_{i3}u(a_{i2}u(a_{i1} u + b_{i1}) + b_{i2})) + b_{i3}\cdots) + b_{iM} 
\]
is the approximation of $F(t, x, u)$ on $U_i$ in the sense that 
\[
|F(t, x, u) - p_i(u)| < C\delta\eps, \quad (t, x) \in U_i, \ \ |u|< \eps. 
\]
The reconstruction of $c(x)$ on $U_i$ is $c_{ii}$ and by the proof of Theorem \ref{thmmain}, we have 
\[
|c(x) - c_{ii}| \leq C \delta \text{ on } U_i.  
\]


\end{document}